\documentclass[12pt]{article}
\usepackage{amsmath, mathrsfs, amsfonts, dsfont, hyperref, amsthm, xspace, cases}

\newtheorem{theorem}{Theorem}[section]
\newtheorem{prop}[theorem]{Proposition}
\newtheorem{remark}[theorem]{Remark}
\newtheorem{definition}[theorem]{Definition}

\newcommand{\msee}{multivalued stochastic evolution equation\xspace}

\newcommand{\Gr}{\ensuremath{\mathrm{Gr}}\xspace}

\def\<{\langle}    \def\>{\rangle}

\renewcommand{\H}{\ensuremath{\mathds{H}}\xspace}
\newcommand{\V}{\ensuremath{\mathds{V}}\xspace}
\newcommand{\VStar}{\ensuremath{\mathds{V^{*}}}\xspace}

\providecommand{\DualInnerProdV}[2]{{}_{\V}\<#1, #2\>_{\V^{*}}}

\providecommand{\ClosedDom}[1][A]{\ensuremath{\overline{D(#1)}}\xspace}

\newcommand{\Zt}{X_{t}-Y_{t}}
\newcommand{\ZtH}{|X_{t}-Y_{t}|_{\H}}
\newcommand{\ZtV}{|X_{t}-Y_{t}|_{\V}}
\newcommand{\ZtS}{|X_{t}-Y_{t}|_{\sigma_{t}}}
\newcommand{\eDKt}{\e^{-\frac{\delta}{2} \omega t}}
\newcommand{\eDKtD}{\e^{-\delta \omega t}}
\newcommand{\zH}{|x-y|_{\H}}
\newcommand{\SpaceNorm}[2]{|#1|_{#2}}
\DeclareMathOperator{\e}{e}

\providecommand{\UnitVect}[1]{\frac{#1}{|#1|}}
\def\E{\mathds E}
\def\P{\mathds P}
\def\Q{\mathds Q}
\def\R{\mathds R}

\def\B{\mathcal B}
\def\C{\mathscr C}
\def\BB{\mathscr B}
\def\F{\mathscr F}
\def\M{\mathcal M}

\hypersetup{
    colorlinks = true,
    linkcolor = black,
    anchorcolor = black,
    citecolor = black,
    filecolor = black,
    pagecolor = black,
    urlcolor = black}

\DeclareFontFamily{U}{bbold}{}
\DeclareFontShape{U}{bbold}{m}{n}{<-5>bbold5<6>bbold6<7>bbold7<8>bbold8
<9>bbold9<9-11>bbold10<11-14>bbold12<14->bbold12}{}
\DeclareSymbolFont{bbold}{U}{bbold}{m}{n}
\DeclareSymbolFontAlphabet{\mathbbold}{bbold}
\DeclareMathSymbol{\Eins}{\mathord}{bbold}{`1}

\title{\Large\textbf{Harnack Inequalities and Applications for Multivalued Stochastic Evolution Equations}}

\author{Shun-Xiang Ouyang\footnote{souyang@math.uni-bielefeld.de; Tel.: +49 521 1064760; fax: +49 521 1066455;}\\
\footnotesize{Department of Mathematics, Bielefeld University,  D-33501 Bielefeld, Germany}
\\
\footnotesize{School of Math. Sci. \& Lab. Math. Com. Sys.,Beijing Normal University, Beijing 100875, China}
}

\begin{document}

\maketitle

\begin{abstract}
By the method of coupling and Girsanov transformation, 
 Harnack inequalities [F.-Y. Wang, 1997] and strong Feller property  are proved for the transition semigroup associated with the multivalued stochastic evolution equation on a Gelfand triple.
The concentration property of the invariant measure for the semigroup is investigated. 
As applications of Harnack inequalities, explicit upper bounds
of the $L^p$-norm of the density,  contractivity, compactness
and entropy-cost inequality for the semigroup are also presented. 
\end{abstract}



\section{Introduction}
Dimension-free Harnack inequality  has been extensively studied (see \cite{Ouyang09, Wang06}).
It is first introduced by Wang \cite{Wang97} 
 for diffusions on Riemannian manifolds with curvature bounded below, then it is considered by
 Aida and Kawabi \cite{AK01, Kaw04, Kaw05}  for some infinite dimensional diffusion processes; 
R\"ockner and Wang \cite{RW03} for generalized Mehler semigroup;
Arnaudon et al. \cite{ATW06} for diffusions on Riemannian manifolds with curvature unbounded below; 
Wang \cite{Wang07} for stochastic porous media equations; 
Da Prato et al.  \cite{DRW09} for singular stochastic equations on Hilbert spaces; 
and the author et al.  \cite{Ouyang09c} for Ornstein-Uhnelbeck processes with jumps in infinite dimensional spaces etc..

For the applications of Harnack inequalities, we refer to 
\cite{BGL01, RW03,RW03b,Wang99,Wang01} 
for contractivity properties and functional inequalities;
\cite{AK01,AZ02, Kaw05} for short time 
heat kernel estimates of infinite dimensional diffusions;
\cite{DRW09} for regularizing properties;
and \cite{BLQ97, GW01} for heat kernel estimates etc..  

In \cite{AK01, Kaw04,Kaw05, RW03, Wang97} etc. 
Harnack inequalities were  proved via gradient estimate. 
In \cite{ATW06} 
the method of coupling and Giranov transformation (\cite{ATW06}) was introduced. 
This method has been applied in  \cite{DRW09, Ouyang09, Ouyang09c, Wang07} etc.. 
In this paper, we adopt this method to establish Harnack inequality for multivalued stochastic 
evolution equations which are generalization of the single valued and  
multivalued linear case considered respectively in \cite{Wang07}.

Multivalued equations attracted the interest of many researchers recently. See
Kr\'ee \cite{Kre82}, C\'epa \cite{Cep94, Cep95, Cep98}, 
Bensoussan and Rascanu \cite{BR97}, 
C\'epa and
L\'epingle \cite{CL97}, and Zhang \cite{Zha07} etc. and references therein.
We prove Harnack inequalities (see Theorem \ref{Thm:MSEE}) and 
strong Feller property (see Theorem \ref{Thm:StrongFeller}) of the invariant measure 
for the transition semigroup associated with the 
multivalued stochastic evolution equations in Banach spaces in the framework of 
\cite{Zha07}. 

As applications of the Harnack inequalities, we prove the invariant measure 
is fully supported on the domain of the underlying multivalued maximal monotone 
operator; and we study the hyperboundedness, ultraboundedness and 
compactness for the transition semigroup. 
We also get a log-Harnack inequality and an entropy-cost inequality. 
We refer to Theorem \ref{Thm:Application} and Theorem \ref{entropy_cost:logHI:MSEE} for details. 

 Zhang \cite[Theorem 5.8]{Zha07} has proved finiteness of the second moment  of the invariant measure of the transition semigroup associated with evolution equations.
 We obtain stronger concentration properties of the invariant measure. 
See Theorem \ref{Thm:Concentration:MSEE}.

The organization of this paper follows. 
We introduce the multivalued stochastic evolution  
equation in Section \ref{Sec:MSEE}; 
study the concentration property of the invariant measure in Section \ref{Sec:Concentration:MSEE}; 
and prove Harnack inequalities in Section \ref{Subsec:MSEE:HI}; then present the applications in Section \ref{Subsec:Application:MSEE}.


\section{Multivalued stochastic evolution equations}
\label{Sec:MSEE}

We first recall the definition of multivalued maximal monotone operator. 
See for instance Br{\'e}zis \cite{Bre73}
for more details.

Denote by $2^{\H}$ for  the set of all subsets of $\H$. 
Let $A\colon \H\to 2^{\H}$ be a 
set-valued operator. 
Define the domain of $A$ by 
\[
       D(A)=\{x\in\H\colon Ax\neq \emptyset\}.
\]
The multivalued operator $A$ is characterized by its {graph}
defined by 
\[      
       \Gr(A)=\{(x,y)\in\H\times\H \colon x\in\H, y\in Ax\}.
\]

\begin{definition}\begin{enumerate}
\item A multivalued operator $A$ on $\H$ is called {monotone} 
if
\[\<x_{1}-y_{1}, x_{2}-y_{2}\>\geq 0, \quad \mathrm{for \ all}\
(x_{1},y_{1}), (x_{2},y_{2})\in\Gr(A).
\]
\item A monotone operator $A$ is called 
	maximal monotone if       
       it must be $(x_{1},y_{1})\in\Gr(A)$
        for any
       $(x_{1},y_{1})\in\H\times \H$ satisfying
       \[
               \<x_{1}-x_{2}, y_{1}-y_{2}\>\geq 0 ,\quad
               \mathrm{for \ all}\ (x_{2},y_{2})\in\Gr(A).
       \]
       That is, $A$ is maximal monotone if its graph $\Gr(A)$ is not contained in
       the graph of any other monotone operator. 
\end{enumerate}
\end{definition}


Let $\V$ be a separable and reflexive Banach space which is continuously and densely 
embedded in a separable Hilbert space $\H$. Then we have an evolution triplet 
$(\V, \H, \VStar)$ satisfying 
\[	\V\subset \H=\H^{*} \subset \VStar	,\]
where $\VStar$ is the dual space of $\V$ and we identify $\H$ with its own dual $\H^{*}$. 

Denote by 
$\SpaceNorm{\cdot}{\V}$, 
$\SpaceNorm{\cdot}{\H}$,
$\SpaceNorm{\cdot}{\VStar}$
the norms in $\V$, $\H$ and $\VStar$ respectively;
by $\<\cdot, \cdot\>_{\H}$ the inner product in \H, and
\(\DualInnerProdV{\cdot}{\cdot}\) the dual relation between $\V$ and $\VStar$. 
In particular, if $v\in\V$ and $h\in\H$, then
\[\DualInnerProdV{v}{h}=\<v, h\>_{\H}.\]


Let $A$ be a multivalued maximal monotone operator on $\H$.
We introduce two sets for every $T>0$:
\begin{enumerate}
\item $\mathscr{V}_{T}(\H)$: the set of all 
$\H$-valued functions of finite variation on $[0,T]$. 
\item $\mathscr{A}_{T}$: the space of all $[u, K]$ such that 
$u\in C([0,T]; \overline{D(A)})$, 
$K\in \mathscr{V}_{T}(\H)$ with $K(0)=0$, 
and for all $x,y\in C([0,T],\H)$ satisfying 
$(x(t),y(t))\in \Gr(A)$, the measure 
$$\<u(t)- x(t), dK(t)-y(t)\,dt\>_{\H}\geq 0.$$
\end{enumerate}

Let $W_{t}$ be a cylindrical Wiener process on $\H$ with respect to 
 a filtered probability space
$(\Omega, \F,  (\F_{t})_{t\geq 0},\P)$. 
Let $B$  
be a single valued operator from 
$\V$ to $\VStar$; and $\sigma$ an operator from $\R_{+}\times \Omega\times 
\H$ to $\H\otimes \H$.
Consider the following \msee
\begin{equation}\label{msee:intro}
\left\{
	\begin{aligned}
		&dX_{t}\in -AX_{t}\,dt  +BX_{t}\, dt+\sigma(t,X_{t})\, dW_{t},\\  
		&X_{0}=x\in \overline{D(A)}.
	\end{aligned}	
\right.
\end{equation}

\begin{definition}\label{Def:Sol:MSSEE}
A pair of $\F_{t}$-adapted random processes $(X_{t}, K_{t})$ is called a solution of Equation 
\eqref{msee:intro} if 
\begin{enumerate}
	\item \([X(\cdot, \omega), K(\cdot, \omega)]\in \mathscr{A}_{T}\) for almost all $\omega \in \Omega$;
	\item For some $q>1$, 
		\(X(\cdot, \omega)\in L^{q}([0,T]; \V)\) for almost all $\omega \in \Omega$;
	\item It holds  for all $t\in[0,T]$ almost surely
		\[X_{t}=X_{0}-K_{t}+\int_{0}^{t} BX_{s}\,ds 
			+\int_{0}^{t}\sigma(s,X_{s})\, dW_{s}.\]		
\end{enumerate}
\end{definition}

The following theorem 
on the existence and uniqueness of the equation \eqref{msee:intro}
 is due to Zhang \cite[Theorem 4.6]{Zha07}.
\begin{theorem}\label{MSEE:Solution:Existence}
Assume the following conditions. 
\begin{description}
	\item[(H1)] \hypertarget{H1:MSEE}
			{$0\in D(A)^{o}$, where $D(A)^{o}$ denotes the interior of $D(A)$;}
	\item[(H2)]  \hypertarget{H2:MSEE}
			{$B$ is hemicontinuous: for every $x,y,z\in\V$, 
			\[[0,1]\ni \varepsilon \mapsto \DualInnerProdV{x}{B(y+\varepsilon z)}\ \mathrm{is\ continuous;}\]}
	\item[(H3)]   \hypertarget{H3:MSEE}
			{For every $x,y\in\V$, 
			\[\DualInnerProdV{x-y}{Bx-By}\leq 0;\]}		
	\item[(H4)]   \hypertarget{H4:MSEE}
			{There exist $\gamma>0$, $\omega\in\R$ and $q>1$ such that for every $x,y\in\V$, 
		\begin{equation}\label{Equation:H4}
			\DualInnerProdV{x-y}{Bx-By}\leq
			  - \gamma |x-y|_{\V}^{q} + \omega |x-y|_{\H}^{2};
		\end{equation}}
	\item[(H5)]   \hypertarget{H5:MSEE}
			{There exists a $C>0$ such that for every $x\in \V$, 
		\[|Bx|_{\VStar}\leq C (1+|x|_{\V}^{q-1}),\]
		where $q$ is the same as in \eqref{Equation:H4};}	
	\item[(H6)]   \hypertarget{H6:MSEE}
		 {Let $\M$ be the set of all progressively measurable sets  with respect to $\F_{t}$. 
		Assume $\sigma$ is $\M\times \B(\H)/\B(\H \otimes \H)$ measurable and 
		there exists a positive constant $C_{\sigma}$ such that for all 
		$(t, \omega)\in\R_{+}\times \Omega$	and $x,y\in\H$, 
		\[
			\begin{aligned}
				&\| \sigma(t,\omega, x) - \sigma(t,\omega, y) \|_{\H\otimes \H} 
						\leq C_{\sigma}|x-y|_{\H},\\
				&\| \sigma(t,\omega, x) \|_{\H\otimes \H} \leq C_{\sigma}(1+|x|_{\H}).		
			\end{aligned}		
		\]}		 
\end{description}
Then there exists a unique solution 
to Equation \eqref{msee:intro} in the sense of Definition \ref{Def:Sol:MSSEE}.
\end{theorem}

The following proposition will play an important role. 
See Zhang \cite[Proposition 3.3]{Zha07} for a proof.
\begin{prop}\label{HI:MSEE:CrucialProp}
Let $[u,K], [\tilde{u},\widetilde{K}]\in \mathscr{A}_{T}$. 
Then the measure
\[
\<u(t)- \tilde{u}(t), dK(t)-d\widetilde{K}(t)\>_{\H}\geq 0,\quad t\in[0,T]
.\]
\end{prop}

\section{Concentration of invariant measures}
\label{Sec:Concentration:MSEE}
Suppose that Conditions 
\hyperlink{H1:MSEE}{(H1)}--\hyperlink{H6:MSEE}{(H6)} hold.
By Theorem \ref{MSEE:Solution:Existence}, Equation 
\eqref{msee:intro} has a unique solution $X_{t}$. 
Define
$$P_{t}f(x)=\E_{\P}f(X_{t})$$
 for every $t\in [0,T]$, $f\in \BB_{b}(\overline{D(A)})$. 
Let $\sigma$ be deterministic and time independent.
Then $P_{t}$ is a Markov semigroup (see \cite[Theorem 5.5]{Zha07}).
 
Zhang \cite[Theorem 5.8]{Zha07} has studied the 
the existence, uniqueness of the invariant measure associated with $P_{t}$. 
He also proved that the invariant measure $\mu$ satisfies 
$\mu(|x|_{\H}^{2})<\infty.$
We prove a stronger concentration properties. 

\begin{theorem}\label{Thm:Concentration:MSEE}
		Assume that \hyperlink{H1:MSEE}{(H1)}--\hyperlink{H6:MSEE}{(H6)} hold
		 with $q\geq 2$. Let $\sigma$ be deterministic and independent of time. 
		Assume further that $\V$ is compactly embedded in $\H$. 
		If $q=2$, then suppose in addition that 
			 $\sigma$ is uniformly bounded and  $\lambda \omega<\gamma$,
	where $\lambda$ is the constant such that $|\cdot|_{\H}\leq \lambda |\cdot|_{\V}$.
	Then there exist an
	invariant measure associated with $P_{t}$ in the sense that
	\[ \int_{\overline{D(A)}} P_{t}f(x)\,\mu(dx) 
		=\int_{\overline{D(A)}} f(x)\,\mu(dx),\quad f\in\BB_{b}(\overline{D(A)}).\]
	Moreover, 
	\begin{equation}\label{Equ:Concentration:MSEE:Moment}
		\int_{\overline{D(A)}} 
			|x|_{\V}^{q}\, \mu(dx)	
				<\infty.		
	\end{equation}		
	If  $\sigma$ is uniformly bounded, then
	for every $q\geq 2$, 
	there exist some $\theta>0$ such that 
	\begin{equation}\label{Equ:Concentration:MSEE:Exp}
		\int_{\overline{D(A)}} 
			\e^{\theta |x|_{\H}^{q}} \mu(dx)	
				<\infty.		
	\end{equation}		   	
\end{theorem}

\begin{proof}
	(1) 	The existence of the invariant measures has been proved 
		in \cite[Theorem 5.8 (i)]{Zha07} for the case $q>2$.
		The extension to the case $q=2$ is not hard. We skip the proof here since the main technique
		can be found below.
		
	(2) From (H3) we know for all $x\in\V$, 
		\begin{equation}\label{Equ:Conc:1}
			\DualInnerProdV{x}{Bx}\leq -\gamma  |x|_{\V}^{q} + \omega|x|_{\H}^{2} +
			\DualInnerProdV{x}{B0}.	
		\end{equation}	
		
		If $q=2$, then
			\begin{equation}\label{Equ:qEqua2:Control}	
				\omega|x|_{\H}^{2}\leq \lambda \omega |x|_{\V}^{2}
					<\gamma  |x|_{\V}^{2}.
			\end{equation}
		If $q>2$, then by Young's inequality, 
		\begin{equation}\label{Equ:qGreater2:Control}
			\omega|x|_{\H}^{2}\leq 
			\lambda \omega |x|_{\V}^{2}
			\leq \frac{2 \varepsilon ^{q}}{q} |x|_{\V}^{q} 
				+ \frac{(\lambda \omega)^{p'}}	{p'\varepsilon^{p'}},
		\end{equation}
		hold for every $\varepsilon>0$, where $p'$ satisfying $1/p'+2/q=1$.
		
		Use the estimate \eqref{Equ:qEqua2:Control} and \eqref{Equ:qGreater2:Control}
		in \eqref{Equ:Conc:1} for $q=2$ and $q>2$ 
		(by taking $\varepsilon$ small enough in this case) respectively, we know
		there are constants $C_{1}, \gamma'>0$ such that 
		\begin{equation}\label{Equ:Conc:2}
			\DualInnerProdV{x}{Bx}\leq C_{1} -2\gamma' |x|_{\V}^{q} +
			\DualInnerProdV{x}{B0}.	
		\end{equation}
		
		By Young's inequality again, we know for any $\tilde{\varepsilon}>0$, 
		\begin{equation}\label{Equ:Conc:3}
			\DualInnerProdV{x}{B0}\leq |x|_{\V}\cdot|B0|_{\VStar}
			\leq \frac{\tilde{\varepsilon}^{q}}{q} |x|^{q}_{\V}
			+\frac{1}{p\tilde{\varepsilon}^{p}}|B0|^{p}_{\VStar},
		\end{equation}
		where $p=\frac{q}{q-1}$. 
		
		Therefore,  we deduce 
		from \eqref{Equ:Conc:2} and \eqref{Equ:Conc:3} by taking $\tilde{\varepsilon}$ small enough to get
		\begin{equation}\label{Equ:Conc:4}
			\DualInnerProdV{x}{Bx}\leq C_{2} -\gamma' |x|_{\V}^{q}	
		\end{equation}
		for some constant $C_{2},\gamma'>0$.
		
		Now we fix a $y\in A0$. Let $(X_{t},Y_{t})$ be the solution to the multivalued 
		stochastic evolution equation \eqref{msee:intro}.
		By definition, we have 
		\begin{equation}\label{XK:control1}
			\<X_t-0, dK(t)-y\,dt\>\geq 0,			
		\end{equation}

		By It\^o's formula, 
		 using  \eqref{Equ:Conc:4} and  \eqref{XK:control1} 
		 and Young's inequality again, we  obtain
		\begin{equation}\label{Equ:Conc:5}
		\begin{aligned}
			& \frac12 d|X_{t}|_{\H}^{2} \\
		    \leq & {}-\DualInnerProdV{X_{t}}{BX_{t}} \,dt 
		    		 -\<X_{t}, dK_{t}\>_{\H}\,dt
			+ \frac12 \|\sigma\|^{2}_{\H\otimes\H}\,dt
		    	+\<X_{t}, \sigma dW_{t}\>\\
		   \leq & \left(C_{3} -\gamma' |X_t|_{\V}^{q}\right)\,dt+|y|_{\H}\cdot |X_{t}|_{\H}\,dt
		   	+ \<X_{t}, \sigma dW_{t}\>\\
		   \leq& \left(C_{4} -\frac{\gamma'}{2} |X_t|_{\V}^{q}\right)\,dt
		      +	 \<X_{t}, \sigma dW_{t}\>,	
		\end{aligned}	
		\end{equation}	
	where $C_{3},C_{4}>0$ are some constants.	
	
	In the calculation of  \eqref{Equ:Conc:5} we also used Young's inequality to 
	get control from $\|\sigma(x)\|_{\H\times\H}\leq C_{\sigma}(1+|x|_{\H})$ 
	if $q$ is strictly greater than 2. 
	If $q=2$, we use the assumption that $\sigma$ uniformly bounded. 

	Therefore, by \eqref{Equ:Conc:5}, we get
	\begin{equation}\label{integrability:int011}
		\int_{0}^{1} \frac{\gamma'}{2} \E^{x} |X_{s}|_{\V}^{q} \,ds
		\leq C_{4} +\frac{1}{2}\left(|x|_{\H}^{2}-\E^{x}|X_{1}|_{\H}^{2}\right).
	\end{equation}
	Consequently we have
	\begin{equation*}
		\int_{0}^{1} P_{s}|\cdot|_{\V}^{q}(x)\,ds  \leq
			\frac{1}{\gamma'} (2C_{4}+ |x|_{\H}^{2}). 
	\end{equation*}
	Hence we obtain 
		$ \mu(|\cdot|_{\V}^{q})<\infty. $
This proves \eqref{Equ:Concentration:MSEE:Moment}.
		
	(3) 	
	For every $\theta>0$, by \eqref{Equ:Conc:5} we  have 
	\begin{equation}\label{Equ:Control:Exp:ItoFormula}
	\begin{aligned} 
			&d\e^{\theta |X_{t}|_{\H}^{q}}\\ 
		=&	\frac{1}{2}\theta q |X_{t}|_{\H}^{q-2}
				  \e^{\theta |X_{t}|_{\H}^{q}} \,d|X_{t}|^{2}_{\H} 
			\\
		&\quad	
			+\frac12 \left( \frac{1}{2}\theta q \e^{\theta |X_{t}|_{\H}^{q}}\right) \left(
					\frac{1}{2}\theta q |X_{t}|_{\H}^{2(q-2)} +\frac{q-2}{2} |X_{t}|_{\H}^{q-4}
				\right)  d\<|X_{t}|_{\H}^{2}, |X_{t}|_{\H}^{2}\>\\
		=& \frac{1}{2}\theta q |X_{t}|_{\H}^{q-2}
				  \e^{\theta |X_{t}|_{\H}^{q}} 
				  \left( d|X_{t}|^{2}_{\H} 	+		 
				  	  2 \theta q |\sigma|_{\H\otimes\H}^{2}
					    |X_{t}|^{q}_{\H} \,dt
					    +(q-2) |\sigma|_{\H\otimes\H}^{2} \,dt	     	
				  \right)\\
		\leq&
			\frac{1}{2}\theta q |X_{t}|_{\H}^{q-2}
			 \e^{\theta |X_{t}|_{\H}^{q}}
		\left(C_{5} -\gamma' |X_t|_{\V}^{q}
			 + 2 \theta q |\sigma|_{\H\otimes\H}^{2}  |X_{t}|^{q}_{\H} \right)\,dt		
			 		      +dM_{t}
	\end{aligned}
	\end{equation}
 for some constant $C_{5}>0$ and some local martingale $M_{t}$. 
 
Since $|\cdot|_{\H}\leq \lambda |\cdot|_{\V}$, for small enough $\theta$, we have 
\begin{equation}\label{Equ:Conc:6}
	d\e^{\theta |X_{t}|_{\H}^{q}} \leq 
	\frac{1}{2}\theta q |X_{t}|_{\H}^{q-2}
			 \e^{\theta |X_{t}|_{\H}^{q}}
		\left(C_{5} -\frac{\gamma'}{2} |X_t|_{\V}^{q}
			  \right)\,dt				      +dM_{t}.
\end{equation} 	

Let us focus at the drift of the right hand of 
\eqref{Equ:Conc:6}. 
By the fact $|\cdot|_{\H}\leq \lambda |\cdot|_{\V}$ and Young's inequality, 
\begin{equation}\label{Equ:Conc:6:drift:control:temp1}
\begin{aligned}
	 \frac{1}{2}\theta q |X_{t}|_{\H}^{q-2}
		\left(C_{5} -\frac{\gamma'}{2} |X_t|_{\V}^{q} \right)
\leq& 	\frac{1}{2}\theta q C_{5} |X_{t}|_{\H}^{q-2}
	 - \frac{1}{2}\theta q \cdot \frac{\gamma'}{2}\lambda^{-q} 
			|X_t|_{\H}^{q} \cdot |X_{t}|_{\H}^{q-2}\\
\leq& C_{6}-\gamma''|X_{t}|_{\H}^{2(q-1)}			
\end{aligned}
\end{equation}
for some constant $C_{6}, \gamma''>0$.

Now let 
\[
	G=\left\{	 |X_{t}|_{\H}^{2(q-1)}\geq 1+\frac{C_{6}}{\gamma''}	\right\}.\]
	Note that on $G^{c}$, both 
	$|X_{t}|_{\H}^{2(q-1)}$ and  $\e^{\theta |X_{t}|_{\H}^{q}}$
	are bounded. 
Therefore 
\begin{equation}\label{Equ:Conc:6:drift:control:temp2}
\begin{aligned}
	&
	\left(C_{6}-\gamma''|X_{t}|_{\H}^{2(q-1)}\right)
	 \e^{\theta |X_{t}|_{\H}^{q}}
			\\
=& -\gamma''{}\left( |X_{t}|_{\H}^{2(q-1)}-\frac{C_{6}}{\gamma''{}}										\right) \e^{\theta |X_{t}|_{\H}^{q}}\\
\leq & -\gamma''{}  \e^{\theta |X_{t}|_{\H}^{q}} \Eins_{G}
		-\gamma''{}\left( |X_{t}|_{\H}^{2(q-1)}-\frac{C_{6}}{\gamma''{}}										\right) 
		\e^{\theta |X_{t}|_{\H}^{q}}
		 \Eins_{G^{c}}\\
\leq & -\gamma''{}  \e^{\theta |X_{t}|_{\H}^{q}} 
	+\gamma''{}  \e^{\theta |X_{t}|_{\H}^{q}} 
		\Eins_{G^{c}}
		-\gamma''{}\left( |X_{t}|_{\H}^{2(q-1)}-\frac{C_{6}}{\gamma''{}}										\right) 
		\e^{\theta |X_{t}|_{\H}^{q}}
		 \Eins_{G^{c}}\\
\leq & C_{7} -\gamma'' 
					 \e^{\theta |X_{t}|_{\H}^{q}}		 		
\end{aligned}
\end{equation}
for some constant $C_{7}>0$. 

Hence from \eqref{Equ:Conc:6:drift:control:temp1} and \eqref{Equ:Conc:6:drift:control:temp2}, 
we can get an estimate of the drift of the right hand side of \eqref{Equ:Conc:6}.
That is, 
\begin{equation}\label{Equ:Conc:7}
	d\e^{\theta |X_{t}|_{\H}^{q}} \leq 
		\left(C_{7} -\gamma'' 
					 \e^{\theta |X_{t}|_{\H}^{q}}
			  \right)\,dt				      +dM_{t}
\end{equation} 

By integrating the inequality \eqref{Equ:Conc:7} from 0 to $n$, we get 
\begin{equation}\label{Equ:Conc:71}
	\e^{\theta |X_{n}|_{\H}^{q}} \leq 
		 \e^{\theta |X_{0}|_{\H}^{q}} +
			C_{7} n-\gamma'' \int_{0}^{n} 
					 \e^{\theta |X_{s}|_{\H}^{q}}\,ds
							      +M_{n}.
\end{equation} 
Then we take expectation for both side of \eqref{Equ:Conc:71} with respect to $\P^{0}$, we get
\begin{equation}\label{Equ:Conc:72}
	\E \e^{\theta |X_{n}|_{\H}^{q}} \leq  1+
		C_{7} n-\gamma'' \int_{0}^{n} 
					\delta_{0}P_{s} \e^{\theta |\cdot|_{\H}^{q}}\,ds.
\end{equation}
It follows that
	\begin{equation}\label{Equ:Conc:8}
		\mu_{n}(\e^{\theta |\cdot|^{q}_{\H}}) 
	\leq \frac{C_{7}}{\gamma''}+\frac{1}{n\gamma''},\quad n\geq 1,
\end{equation}
where 
\[\mu_{n}=\frac{1}{n}\int_{0}^{n} \delta_{0}P_{s}\,ds,\quad n\geq 1.\]
Note that $\mu$ is the weak limit of $\mu_{n}$ (refer to the proof of \cite[5.8]{Zha07}),
we can deduce from \eqref{Equ:Conc:8} to get 
$ \mu(\e^{\theta |\cdot|_{\H}^{q}})<\infty.$
This proves \eqref{Equ:Concentration:MSEE:Exp}.
\end{proof}


\section{Harnack inequalities}
\label{Subsec:MSEE:HI}
In the following we assume conditions \hyperlink{H1:MSEE}{(H1)}--\hyperlink{H5:MSEE}{(H5)}
 in Theorem \ref{MSEE:Solution:Existence} and instead of 
\hyperlink{H6:MSEE}{(H6)} we suppose that 

\begin{description}
\item[(H6$'$)] 
\hypertarget{H6prime:MSEE}
{$\sigma\colon [0,\infty)\times \Omega\to \H\otimes\H$ be a nondegenerate Hilbert-Schmidt operator 
	uniformly bounded in time $t\in[0,\infty)$ and $\omega\in\Omega$.}
\end{description}

For every $x\in\H$,  define
\[ |x|_{\sigma_{t}}=
\begin{cases}
	|y|_{\H} \quad &\mathrm{if}\ x=\sigma_{t}y \ \mathrm{for\ some}\ y\in\H,\\
	\infty,		&\mathrm{otherwise}	. 
\end{cases}
\]
The distance associated with $|\cdot|_{\sigma_{t}}$ is called the intrinsic distance
 induced by $\sigma_{t}$.

 We are going to prove the following Harnack inequality for the semigroup $P_{t}$
 associated with the solution process of Equation \eqref{msee:intro}. 

\begin{theorem}\label{Thm:MSEE}
		Assume \hyperlink{H1:MSEE}{(H1)}--\hyperlink{H5:MSEE}{(H5)} and 
	\hyperlink{H6prime:MSEE}{(H6\,$'$)}.
	Suppose that there exists some nonnegative constant $r\geq q-4$, and
	some strictly positive continuous function $\zeta_{t}$ on $[0,\infty)$
	such that 
	\begin{equation}\label{ControlNormCond:MSEE}
	  \zeta_{t}^{2}|x|_{\sigma_{t}}^{2+r}\cdot |x|_{\H}^{q-2-r}
		\leq |x|_{\V}^{q} ,\quad\mathrm{for\ all}\ x\in\V, \ t\geq 0	
	\end{equation}
holds on $\Omega$. 
Then for every $T>0$, $\alpha>1$, $x,y\in \overline{D(A)}$
 and $f\in\C_{b}^{+}(\overline{D(A)})$,  
the following inequality holds
\begin{equation}\label{HI:MSEE}
(P_{T}f)^{\alpha}(x) \leq 
\exp\left(\frac{\alpha}{2(\alpha-1)} \Theta_{T}|x-y|^{\frac{2(4+r-q)}{2+r}}_{\H}	\right)
P_{T}f^{\alpha}(y),
\end{equation}	
where 
\begin{equation}\label{ThetaTDef}
\Theta_{T}=
4\delta^{-\frac{2(3+r)}{2+r}}\gamma^{-\frac{2}{2+r}}				  
	\frac{\left(  
				\int_{0}^{T} \zeta_{t}^{2}\eDKtD\,dt	
					  \right)^{\frac{r}{2+r}}
					  }
					  {\left(\int_{0}^{T} \zeta_{t}\eDKtD \,dt\right)^{2}}
\end{equation}
with
\begin{equation}\label{HI:MSEE:Delta}
\delta=1-\frac{q}{4+r}.
\end{equation}

Assume the diffusion coefficient $\sigma$ is independent of $(t,\omega)$
and the function $\zeta_{t}$ in  \eqref{ControlNormCond:MSEE} 
is taken as constant $\zeta$. Then  \eqref{HI:MSEE} holds with $\Theta_{T}$ replaced by
\begin{equation}\label{ThetaTDef2}
\widetilde{\Theta}_{T}=
4\delta^{-1}	\gamma^{-\frac{2}{2+r}}			   				  		
		\zeta^{-\frac{4}{2+r}} \left[\omega^{-1}(1-\e^{-\delta \omega T })\right]^{-\frac{4+r}{2+r}}.
\end{equation}
\end{theorem}

\vspace{-0.5em}
\begin{proof}
The proof is divided into six steps. We outline 
the main procedure of the proof of \eqref{HI:MSEE} 
in the first step and then realize the idea in the next four steps. 
The simplification from \eqref{ThetaTDef} to \eqref{ThetaTDef2} 
is obtained in the last step. 

(1) \emph{Main Idea}. 

Consider the following coupled \msee 
\begin{subequations}\label{CoupledEqu:MSEE1}
\begin{numcases}{}
	dX_{t}\in -AX_{t}\,dt + BX_{t}\,dt+\sigma(t)\,dW_{t}-
			 U_{t}\,dt,	                                                        \label{CoupledEqu:MSEE11}	\\
	dY_{t}\in -AY_{t}\,dt + BY_{t}\,dt+	\sigma(t)\,dW_{t}\label{CoupledEqu:MSEE12} 
\end{numcases}
\end{subequations}
with initial conditions \(X_{0}=x\in\overline{D(A)},\ Y_{0}=y\in\overline{D(A)},\) 
and the drift $U_{t}$ in \eqref{CoupledEqu:MSEE11} is of the following form 
\begin{equation}\label{Equation:Drift:HIMSE:Ut}
U_{t}=\frac{\eta_{t}(X_{t}-Y_{t})}{\ZtH^{\delta}}\Eins_{\{t<\tau\}},
\end{equation}
where the stopping time $\tau$ in \eqref{Equation:Drift:HIMSE:Ut} 
is the coupling time of $X_{t}$ and $Y_{t}$ defined by
\[\tau=\inf\{t\geq 0\colon X_{t}=Y_{t}\},\] 
 the power $\delta$ in \eqref{Equation:Drift:HIMSE:Ut}
  is a constant in $(0,1)$ (see \eqref{HI:MSEE:Delta}) and 
$\eta_{t}$ is a deterministic function on $[0,\infty)$.  
Both $\delta$ and $\eta_{t}$ in \eqref{Equation:Drift:HIMSE:Ut} will be specified later such that 
the following two crucial conditions 
			\begin{equation}\label{A1CrucialFact}
					X_{T}=Y_{T}	\quad \textrm{a.s.}
			\end{equation}	
and 
		\begin{equation}\label{A2CrucialFact}
			\E_{\P}\exp \left( \int_{0}^{T} \frac{\eta_{t}^{2}}{2}
			\frac{\ZtS^{2}}{\ZtH^{2\delta}}\Eins_{\{t<\tau\}} \,dt
				 \right) <\infty.
		 \end{equation}
are satisfied.

Let $N_t=\int_{0}^{t} \<\sigma_{s}^{-1}U_{s}, dW_{s}\>.$
By \eqref{A2CrucialFact} we know 
\[R_{t}=\exp\left( N_t
				-\frac12 [N]_t
		     \right)
, \quad t\in[0,T]\]
is a martingale on $(\Omega, \F_{T}, (\F_{t})_{0\leq t\leq T}, \P)$.
Then we can define a new probability measure $\Q$ on 
$(\Omega, \F_{T})$
by setting 
\(\Q|_{\F_{T}}=R_{T}\P.\)
By Girsanov's theorem, 
\[\widetilde{W}_{t}:=W_{t}-\int_{0}^{t} \sigma_{s}^{-1} U_{s}\,ds\]
is still a cylindrical Wiener process on $(\Omega, \F_{T}, (\F_{t})_{0\leq t\leq T}, \Q)$. Hence Equation 
\eqref{CoupledEqu:MSEE11} can be rewritten in the following way
\[
	dX_{t}\in -AX_{t}\,dt + BX_{t}\,dt+\sigma(t)\,d\widetilde{W}_{t}
\]
with initial condition $X_{0}=x$. 

By the uniqueness of the solution, the transition law of 
$(X_{t})_{t\in[0,T]}$ under $\Q$ is the same with the transition law of $(Y_{t})_{t\in[0,T]}$
under $\P$. 
So by the fact \eqref{A1CrucialFact} which will be verified, we have
\begin{equation}\label{SemigroupAnotherFormula}
	P_{T}f(x)=\E_{\Q}f(X_{T})=\E_{\Q}f(Y_{T})=\E_{\P}Rf(Y_{T})	.
\end{equation}
Note that we also have 
\(P_{T}f(y)=\E_{\P}f(Y_{T}),\)
therefore by applying H\"older's inequality to \eqref{SemigroupAnotherFormula}, we  get
\begin{equation}\label{HI:MSEE:HItemp1}
(P_{T}f)^{\alpha} (x)\leq \bigl(\E_{\P}R_{T}^{\alpha/(\alpha-1)}\bigr)^{\alpha-1} P_{T}f^{\alpha}(y).
\end{equation}
Then the proof is completed by estimating $\E_{\P}R_{T}^{\alpha/(\alpha-1)}$.

(2) \emph{Existence of the solution of the coupled equation \eqref{CoupledEqu:MSEE1}}.

Note that the function 
\[(u,v)\mapsto \frac{u-v}{|u-v|_{\H}^{\delta}}\]
satisfies the monotone condition off the diagonal (see \cite[Appendix A]{Wang07}).

By applying Theorem \ref{MSEE:Solution:Existence} we see 
the coupled equation \eqref{CoupledEqu:MSEE1}
has a solution up to the coupling time $\tau$. 
So there exists continuous processes $(X,K)\in\mathscr{A}_{T\wedge \tau}$
and $(Y,\widetilde{K})\in\mathscr{A}_{T\wedge \tau}$ such that 
for all $t<\tau$, 
\begin{subequations}\label{IntegralCoupledEqu:MSEE}
\begin{numcases}{}
	X_{t}=x-K_{t}+\int_{0}^{t} BX_{s}\,ds 
			+\int_{0}^{t}\sigma(s)\, dW_{s}-\int_{0}^{t} U_{s}\,ds,
			\label{IntegralCoupledEqu:MSEEa} \\ 
	Y_{t}=y-\widetilde{K}_{t}+\int_{0}^{t} BY_{s}\,ds +\int_{0}^{t}\sigma(s)\, dW_{s}.		
				\label{IntegralCoupledEqu:MSEEb}
\end{numcases}
\end{subequations}

On the other hand,  it is obvious that the solution of Equation \eqref{CoupledEqu:MSEE12} (or equivalently, Equation \eqref{IntegralCoupledEqu:MSEEb})
can be extended to be a solution for all time $t\in [0,T]$. 
Let $(Y_{t},\widetilde{K})_{t\geq 0}$ solves Equation \eqref{CoupledEqu:MSEE12}. 
Now we get solution of  Equation \eqref{CoupledEqu:MSEE11}
(or \eqref{IntegralCoupledEqu:MSEEa})
 by defining 
\(X_{t}=Y_{t}\), \(K_{t}=\widetilde{K}_{t}\)
for all $t\geq \tau$.

(3) \emph{Verify \eqref{A1CrucialFact}}.

Apply It\^o's formula (see e.g. \cite{KR07}) or  Zhang \cite[Theorem A.1]{Zha07} etc.) to 
\(\sqrt{\ZtH^{2}+\varepsilon}\) and then let $\varepsilon \downarrow 0$, 
by using condition (H4) we have for $t<\tau$
\[
\begin{aligned}
	d\ZtH^{2}\leq &{}-\<\Zt, dK_{t}-d\widetilde{K}_{t}\>_{\H}\,dt\\
			&+(-\gamma \ZtV^{q}+\omega\ZtH^{2}-\eta_{t}\ZtH^{2-\delta})\,dt.
\end{aligned}
\]
By Proposition \ref{HI:MSEE:CrucialProp}, for all $t<\tau$ we have 
\[
d\ZtH^{2}\leq 
			(-\gamma \ZtV^{q}+\omega\ZtH^{2}-\eta_{t}\ZtH^{2-\delta})\,dt.
\]
Then 
\begin{equation}\label{HI:MSEE:DerivativeNeed2}
d\left( \ZtH^{2}\e^{-\omega t} \right)\leq -\e^{-\omega t}
	 \left( \gamma \ZtV^{q} +\eta_{t} \ZtH^{2-\delta} \right)\,dt.
\end{equation}
Hence by \eqref{HI:MSEE:DerivativeNeed2} we get
\begin{equation}
\label{Equ:MSE:FinalInequality:StoppingTime}
\begin{aligned}
	 d\left( \ZtH^{2} \e^{-\omega t} \right)^{\delta/2}
   \leq &\,\frac{\delta}{2} \left( \ZtH^{2} \e^{-\omega t} \right)^{\delta/2-1}
	   	\cdot\left( -\e^{-\omega t} \eta_{t} \ZtH^{2-\delta} \right)\,dt	\\
    =& -\frac{\delta}{2}\eDKt\eta_{t}\,dt. 		
\end{aligned}
\end{equation}

We take 
\begin{equation}\label{MSEHI:eta_t:specify}
\eta_{t} =\vartheta_{T} \zeta_{t} \eDKt
\end{equation}
with 
\[
	\vartheta_{T}
	=\frac{2\delta^{-1}\zH^{\delta}} {\int_{0}^{T} \zeta_{t}\eDKtD \,dt}.
\]
It must be $T\geq \tau$ and hence $X_{T}=Y_{T}$.  
Otherwise, if $T<\tau$,  by taking integral from 0 to $T$ for 
both sides of the inequality 
\eqref{Equ:MSE:FinalInequality:StoppingTime}, we can obtain
\begin{equation}\label{XYmeetAtTtemp1}
	|X_{T}-Y_{T}|_{\H}^{\delta}\e^{-\frac{\delta}{2}\omega T} \leq \zH^{\delta} 
		-\frac{\delta}{2}\int_{0}^{T} \eDKt \eta_{t}\,dt.
\end{equation}
By \eqref{MSEHI:eta_t:specify}
the right hand side of \eqref{XYmeetAtTtemp1} equals to 0. 
So we can conclude $X_{T}=Y_{T}$ from \eqref{XYmeetAtTtemp1}. 
But this is contradict with the assumption
that $T<\tau$. 

(4) \emph{Verify \eqref{A2CrucialFact}}.

From \eqref{HI:MSEE:DerivativeNeed2} and the assumption 
\eqref{ControlNormCond:MSEE} we can get for all $t\leq \tau$
\begin{equation}\label{VerifyA2temp1}
\begin{aligned}
       d\left( \ZtH^{2}\e^{-\omega t} \right)^{\delta}
=& \delta 	\left( \ZtH^{2}\e^{-\omega t} \right)^{\delta-1} d\left( \ZtH^{2}\e^{-\omega t} \right) \\
\leq& -\delta \gamma \eDKtD \ZtH^{2(\delta-1)}\cdot \ZtV^{q}\,dt\\
\leq& -\delta \gamma \zeta_{t}^{2} \eDKtD 
 		\frac{\ZtS^{2+r}}{\ZtH^{2+r-2(\delta-1)-q}}\,dt. 
\end{aligned}
\end{equation}

Let $\delta$ be defined as in \eqref{HI:MSEE:Delta}. 
Then from \eqref{VerifyA2temp1} we get
\begin{equation}\label{VerifyA2temp101}
d\left( \ZtH^{2}\e^{-\omega t} \right)^{\delta} \leq 
 -\delta \gamma \zeta_{t}^{2} \eDKtD 
 		\frac{\ZtS^{2+r}}{\ZtH^{\delta(2+r)}}\,dt. 
\end{equation}
According to \eqref{MSEHI:eta_t:specify}, we have 
\(
\zeta_{t}^{2}  =\frac{\eta_{t}^{2}}{\vartheta_{T}^{2}}  \eDKtD
\).
Now by integrating  both sides of the inequality \eqref{VerifyA2temp101} from $0$ to $T$, 
 we get (note that $X_{T}=Y_{T}$)
\[
\frac{\delta \gamma}{\vartheta_{T}^{2}} 
\int_{0}^{T} \frac{\eta_{t}^{2}\ZtS^{2+r}}{\ZtH^{\delta(2+r)}}
\,dt\leq \zH^{2\delta}.
\]
By H\"older's inequality, we have
\begin{equation*}
\begin{aligned}
	 \int_{0}^{T} \frac{\eta_{t}^{2}\ZtS^{2}}{\ZtH^{2\delta}}\,dt
\leq    &\left(   \int_{0}^{T}
				\frac{\eta_{t}^{2}\ZtS^{2+r}}{\ZtH^{\delta(2+r)}}\,dt
					  \right)^{\frac{2}{2+r}}
    		\left(  
				\int_{0}^{T} \eta_{t}^{2}\,dt	
					  \right)^{\frac{r}{2+r}}\\
\leq &\left(
		\frac{\vartheta_{T}^{2}}{\delta \gamma} 
		\zH^{2\delta}
		\right)^{\frac{2}{2+r}}	
		\cdot \vartheta_{T}^{\frac{2r}{2+r}}
		\left(  
				\int_{0}^{T} \zeta_{t}^{2}\eDKtD\,dt	
					  \right)^{\frac{r}{2+r}}
					.
\end{aligned}					  			  		  
\end{equation*}
Note that 
\[
\vartheta_{T}^{2}=\frac{4\delta^{-2}\zH^{2\delta}} {\left(\int_{0}^{T} \zeta_{t}\eDKtD \,dt\right)^{2}},	
\]
we obtain
\begin{equation}\label{VerifyEstimateTemp2}
\begin{aligned}
 	\int_{0}^{T} \frac{\eta_{t}^{2}\ZtS^{2}}{\ZtH^{2\delta}}\,dt 
\leq& 4\delta^{-\frac{2(3+r)}{2+r}}\gamma^{-\frac{2}{2+r}}					  
	\frac{\left(  
				\int_{0}^{T} \zeta_{t}^{2}\eDKtD\,dt	
					  \right)^{\frac{r}{2+r}}
					  }
					  {\left(\int_{0}^{T} \zeta_{t}\eDKtD \,dt\right)^{2}}
		\zH^{\frac{2(4+r-q)}{2+r}}					  
  .
\end{aligned}				  
\end{equation}
Now it is clear  that \eqref{A2CrucialFact} holds. 

(5) \emph{Estimate of\/ \(\E R_{T}^{\alpha/(\alpha-1)}\)}.

Denote $\beta=\alpha/(\alpha-1)$. 
Since $R_{t}$ is a $\P$-martingale, 
for any $p,q>1$ with $q=p/(p-1)$, 
we have
\begin{equation}
\label{HI_temp_cb2}
\begin{split}
&\E_{\P}R_T^{\alpha/(\alpha-1)}\\
=&\E_{\P}\exp\left(\beta N_T-\frac{1}{2}p\beta^2[N]_T\right) 
		\exp\left(\frac{1}{2}\beta(\beta p-1)[N]_T\right)\\
\leq& \left[\E_{\P}\exp\left(p\beta N_T-\frac{1}{2}p^2\beta^2[N]_T\right)\right]^{1/p}\cdot
\left[\E_{\P}\exp\left(\frac{1}{2}\beta q(\beta p-1)[N]_T\right)\right]^{1/q}\\
=&\left[\E_{\P}\exp\left(\frac{1}{2}\beta q(\beta p-1)[N]_T\right)\right]^{1/q}.
\end{split}
\end{equation}
Hence we get \eqref{HI:MSEE} by using the estimate \eqref{VerifyEstimateTemp2} 
and letting $p$ go to 1.

(6) Suppose that $\sigma$ is independent of $(t,\omega)$. And we take $\zeta_{t}$ in  \eqref{ControlNormCond:MSEE} 
as a constant $\zeta$. 
 Then we can simplify $\Theta_{T}$ as follows. 
\begin{equation*}
\begin{aligned}
\Theta_{T}=
& 4\delta^{-\frac{2(3+r)}{2+r}}\gamma^{-\frac{2}{2+r}}		\frac{\left(  \zeta^{2}
				\int_{0}^{T} \eDKtD\,dt
					  \right)^{\frac{r}{2+r}}
					  }
					  {\left( \zeta \int_{0}^{T} \eDKtD \,dt\right)^{2}}\\
= & 4\delta^{-\frac{2(3+r)}{2+r}}\gamma^{-\frac{2}{2+r}}	
		\zeta^{-\frac{4}{2+r}} 					  				  
	\bigl[  (\delta \omega)^{-1}
			(1-\e^{-\delta \omega T })
					  \bigr]^{-\frac{4+r}{2+r}}					  
\\
= & 4\delta^{-1}	\gamma^{-\frac{2}{2+r}}			   				  		\zeta^{-\frac{4}{2+r}} \left[\omega^{-1}(1-\e^{-\delta \omega T })
					\right]^{-\frac{4+r}{2+r}} 
=\widetilde{\Theta}_{T}					
\end{aligned}					  
\end{equation*}

\end{proof}

\begin{remark}
We refer to \cite[Corollary 1.3]{Wang07} for sufficient conditions for 
\eqref{ControlNormCond:MSEE}.
\end{remark}

\begin{remark}
Consider the following multivalued stochastic differential equation
on $\R^d$
\begin{equation*}\label{Equ:MSDE:X}
       dX_{t}+AX_{t}dt \ni b(X_{t})dt+ dW_{t},\quad X_{0}=x\in\overline{D(A)},
\end{equation*}
where $A$ is a maximal monotone operator on $\R^{d}$ with
$D(A)^{o}\neq \emptyset$, $W_{t}$ is a Wiener process on a filtered
probability space $(\Omega, \F, (\F_{t})_{t\geq 0}, \P)$, 
$b\colon \R^{d}\to\R^{d}$ is a 
 continuous function such that there exists some $\omega \in\R$
\begin{equation*}\label{MSDSE:NonLip:Drift}
 \<x-y, b(x)-b(y)\>\leq \omega|x-y|^2, \quad x,y\in\R^d.
\end{equation*} 

Then by \cite[Thorem 2.8]{RWZ08},
the solution $X_{t}$ exists. 
For the  transition semigroup $P_t$ associated with $X_t$, 
we have 
\begin{equation}\label{Equ:HI:MSDE:Monotonic}
	(P_{t}f)^{\alpha}(x)\leq  
				\exp\left(
				\frac{\alpha \omega |x-y|^{2}}{(\alpha-1)(1-\e^{-2\omega t})}
			\right)
				 P_{t}f^{\alpha}(y).
\end{equation}
for every $t>0$, $\alpha>1$, $x,y\in \overline{D(A)}$ and
$f\in\C_{b}^{+}(\overline{D(A)})$.

The proof is similar to the proof of  Theorem \ref{Thm:MSEE} except that we choose another simpler drift. 
That is, instead of \eqref{CoupledEqu:MSEE1}, we consider the following coupled equations
\begin{subequations}\label{HI:MSDE:Couple:Monotonic}
	\begin{numcases}{}
		dX_t+AX_{t}\,dt \ni dW_t+b(X_t)\,dt-
			\xi_t |x-y| \UnitVect{X_t-Y_t} \Eins_{\{t<\tau\}}\,  dt, 
					\label{HI:MSDE:Couple:Monotonic1}\\
		dY_t+AY_{t}\,dt \ni dW_t+b(Y_t)\,dt,  \label{HI:MSDE:Couple:Monotonic2}
	\end{numcases}
\end{subequations}
with initial data $X_{0}=x$ and $Y_{0}=y$, 
where $\tau$ is the coupling time of $X_{t}$ and $Y_{t}$
and 
\begin{equation*}
		\xi_{t}=\frac{\e^{-\omega t}}{\int_{0}^{T}\e^{-2\omega s}\,ds},\quad t\geq 0.
	\end{equation*}
\end{remark}


\section{Applications}\label{Subsec:Application:MSEE}
Zhang \cite[Corollary 5.3]{Zha07} studied 
Feller property of $P_{t}$. We prove strong Feller property
for $P_{t}$ under additional conditions.

\begin{theorem}\label{Thm:StrongFeller}
	Assume \hyperlink{H1:MSEE}{(H1)}--\hyperlink{H5:MSEE}{(H5)}, 
	\hyperlink{H6prime:MSEE}{(H6\,$'$)} and \eqref{ControlNormCond:MSEE}
with $q<4+r$. Then 
the semigroup $P_{t}$ is strongly Feller and 
 the following estimate holds 
\begin{equation}\label{MSEE:StrongFellerEstimate}
|P_{t}f(x)-P_{t}f(y)|\leq 
 \|f\|_{\infty}
\Theta_{t}^{1/2} \zH^{\frac{4+r-q}{2+r}}\cdot
       		\exp\left(\frac{1}{2}\Theta_{t} \zH^{\frac{2(4+r-q)}{2+r}}\right)
\end{equation}
for every $t>0$, $x,y\in\ClosedDom$ and $f\in\BB_{b}(\ClosedDom)$.
Moreover, if the invariant measure $\mu$ exist, then for every $f$ in $L^{p}(\ClosedDom,\mu)$ with $p>1$,  $P_{t} f$ is continuous on
$\ClosedDom$. 
\end{theorem}

\begin{proof}
The last statement follows directly from  the Harnack inequality \eqref{HI:MSEE}
by applying \cite[Proposition 4.1]{DRW09}. 
We only need to show the estimate \eqref{MSEE:StrongFellerEstimate}.

We follow the notation in the proof of Theorem \ref{Thm:MSEE}.
By \eqref{SemigroupAnotherFormula}, we see
\begin{equation}\label{MSEE:StrongFeller1}
\begin{aligned}
          |P_{T}f(x)-P_{T}f(y)|
       &=|\E_{\Q}f(X_{T}) -\E_{\P}f(Y_{T}) | 
       	 =|\E_{\P}R_{T}f(X_{T}) -\E_{\P}f(X_{T}) |  \\ 
       &=\E_{\P}|f(X_{T}) (1-R_{T})| \leq \|f\|_{\infty} \E_{\P}|1-R_{T}|.  
\end{aligned}
\end{equation}

It is clear 
\begin{equation}\label{MSEE:StrongFeller100}
	   \bigl( \E_{\P}| 1-R_{T}| \bigr)^{2} 
       	    \leq  \E_{\P}(1-R_{T})^{2} =\E_{\P}R_{T}^{2}-1.
\end{equation}
By \eqref{HI_temp_cb2}, we know 
\begin{equation}\label{MSEE:StrongFeller101}
\E_{\P} R_{T}^{2}\leq \exp\left(\Theta_{T} \zH^{\frac{2(4+r-q)}{2+r}}\right).
\end{equation} 
Using the elementary inequality 
\(\e^{r}-1\leq r\e^{r}\), $r\geq 0$, 
we can deduce from \eqref{MSEE:StrongFeller100} and \eqref{MSEE:StrongFeller101} to get 
\begin{align*}
	   \bigl( \E_{\P}| 1-R_{T}| \bigr)^{2} 
       \leq & \Theta_{T} \zH^{\frac{2(4+r-q)}{2+r}}
       		\exp\left(\Theta_{T} \zH^{\frac{2(4+r-q)}{2+r}}\right).
\end{align*}
Substitute the estimate above into \eqref{MSEE:StrongFeller1} we can 
obtain \eqref{MSEE:StrongFellerEstimate}.
\end{proof}

\begin{theorem}\label{Thm:Application}
	Assume \hyperlink{H1:MSEE}{(H1)}--\hyperlink{H5:MSEE}{(H5)}, 
	\hyperlink{H6prime:MSEE}{(H6\,$'$)}, \eqref{ControlNormCond:MSEE}, 
	$\sigma$ is independent of $(t,\omega)$ and the invariant measure $\mu$ 
	of the semigroup $P_{t}$ exist.
Then 
\begin{enumerate}
\item The invariant measure $\mu$ is fully supported on $\overline{D(A)}$.
\item For every $x\in\overline{D(A)}$, $t>0$, the
	transition density $p_{t}(x,\cdot)$ (with respect to $\mu$) exist and for
	every $\alpha>1$ 
	\[
		\|p_{t}(x,\cdot)\|_{L^{\alpha}(\overline{D(A)},\mu)} \leq 
		\left[
			\int_{\ClosedDom}\exp\left(-\frac{\alpha}{2}\widetilde{\Theta}_{t} \zH^{\frac{2(4+r-q)}{2+r}} 								\right)\,\mu(dy)
		\right]^{-(\alpha-1)/\alpha}.
	\]
\item Suppose $K\leq 0$. 
	\begin{enumerate}
		\item If $q=2$ and 
			  $\lambda \omega<\gamma$,
			  where $\lambda$ is the constant such that $|\cdot|_{\H}\leq \lambda |\cdot|_{\V}$, 
			then $P_{t}$ is hyperbounded .
		\item If $q>2$, then $P_{t}$ is ultrabounded. More precisely, there 
				exist some constant $c>0$ such that 
			\begin{equation}\label{Equ:Ultrabound:MSEE}
				\|P_{t}\|_{2\rightarrow \infty } \leq \exp\bigl(c(
					1+t^{-\frac{q}{q-2}}
							)\bigr).	
			\end{equation}						
	\end{enumerate} 
	Consequently,  $P_{t}$ is compact for large $t>0$ for both cases. 
\end{enumerate}
\end{theorem}
\begin{proof}
(1) If $\operatorname{supp}\mu \neq \ClosedDom$, then there exists some 
    $x_{0}\in \ClosedDom$ and $r>0$ such that \(\mu(B_{r}(x_{0}))=0\), 
    where $B_{r}(x_{0})=\{y\in\ClosedDom\colon |y-x_{0}|\leq r\}$.
    
    Applying \eqref{HI:MSEE} to the function $\Eins_{B_{r}(x_{0})}$ for $\alpha=2$ and
    $t\geq 0$, we have
    \begin{equation}\label{HIMSEE:Apply:1:Equ0}
    	(P_{t} \Eins_{B_{r}(x_{0})})^{2} (x)
	\exp\left( -\widetilde{\Theta}_{t} \zH^{\frac{2(4+r-q)}{2+r}} \right)
	\leq P_{t}  \Eins_{B_{r}(x_{0})} (y).
    \end{equation}
    Hence,  by integrating both sides of \eqref{HIMSEE:Apply:1:Equ0}
    	 with respect to $\mu(dy)$, we  obtain
      \[
      \begin{aligned}
    	&(P_{t} \Eins_{B_{r}(x_{0})})^{2}(x) 
	\int_{\ClosedDom}
	\exp\left( -\tilde{\Theta}_{t} \zH^{\frac{2(4+r-q)}{2+r}} \right)\mu(dy)\\
	\leq& \mu(P_{t}  \Eins_{B_{r}(x_{0})})=\mu(\Eins_{B_{r}(x_{0})})=0. 
    \end{aligned}
    \]
   This implies $P_{t}(x_{0}, B_{r}(x_{0}))=0$ for all $t\geq 0$. 
   Therefore, 
   \begin{equation}\label{MSEE:APP:FullSupport1}
   \P( | X_{t}(x_{0})-x_{0} |_{\H}\leq r)=0,\quad t>0,
   \end{equation}
   where $X_{t}(x_{0})$ denotes the solution to \eqref{msee:intro}
   with $X_{0}(x_{0})=x_{0}$. 
   
   Since $X_{t}$ is continuous on $\H$, by letting $t\to 0$ in
     \eqref{MSEE:APP:FullSupport1}, we have
    \[
       \P(| X_{0}(x_{0})-x_{0} |_{\H}\leq r)=0.
    \]
   But obviously this is impossible. So it must be $\operatorname{supp}\mu = \ClosedDom$.
   
(2) We refer to \cite{RW03} or \cite{Wang07} for the idea of the proof.

(3) Since $K\leq 0$, for any $t> 0$, we know 
		\[\frac{K}{1-\e^{-\delta Kt}}\leq \frac{1}{\delta t}.\]
      Therefore, by Theorem \ref{Thm:MSEE}, there exist some constant $C_{8}$ 	
      such that for every $ x,y\in\overline{D(A)}$ and $ \ t>0$, 
      \begin{equation}\label{MSEE:BBS:1}
      		(P_{t}f)^{2}(x)\exp
		\left( -\frac{C_{8}|x-y|^{\frac{2(4+r-q)}{2+r}}_{\H}}{t^{\frac{4+r}{2+r}}} \right)
		\leq P_{t}f^{2}(y), 
      \end{equation} 
      where $f\in L^{2}(\overline{D(A)}, \mu)$ with $\mu(f^{2})=1$.
      
      By integrating both sides of \eqref{MSEE:BBS:1}  with respect to  $\mu(dy)$
      over $B_{1}(0)=\{x\in\overline{D(A)}\colon |x|_{\H}\leq 1\}$, we see
      for every $ x\in\overline{D(A)}$ and $ \ t>0$, 
	    \begin{equation}\label{MSEE:BBS:2}
	      		(P_{t}f)^{2}(x)\leq 
			\frac{1}{\mu\bigl(B_{1}(0)\bigr)}
			\exp
			\left( \frac{C_{8}\bigl(1+|x|_{\H}\bigr)^{\frac{2(4+r-q)}{2+r}}}{t^{\frac{4+r}{2+r}}} \right).
	\end{equation}
	
	(i)	 
	If $q=2$, then by taking square and integration with respect to $\mu(dx)$ for
	  both sides of the equation \eqref{MSEE:BBS:2}, and using 
	 Theorem \ref{Thm:Concentration:MSEE}, we have 
	 \[
	 	\int_{\overline{D(A)}} (P_{t}f)^{4}(x)\,\mu(dx)\leq \frac{1}{\mu\bigl(B_{1}(0)\bigr)}
				\int_{\H} \exp\left(  
						\frac{C_{8}\bigl(1+|x|_{\H}\bigr)^{2}}{t^{\frac{4+r}{2+r}}}
				  \right)\,\mu(dx)<\infty
	 \]
	for $t>0$ big enough. This proves 
	\(\|P_{t}\|_{2\rightarrow 4} <\infty \)
	for sufficiently big $t>0$. That is, $P_{t}$ is hyperbounded. 
	
	(ii)
	Assume $q>2$.  Then the inequality \eqref{Equ:Conc:6}
	implies 
	\begin{equation}\label{Equ:Conc:61}
	d\e^{\theta |X_{t}|_{\H}^{q}} \leq 
		\left(C_{9} -{\gamma'''} |X_t|_{\H}^{2(q-1)} 		
		 \e^{\theta |X_{t}|_{\H}^{q}}
			  \right)\,dt				      +dM_{t}
	\end{equation} 
	for some constant $C_{9}, \theta, \gamma'''>0$.
	
	Let $g(t)$ be the solution to the following equation
	\[
		dg(t)=\left(		C_{9}-\gamma''' \theta^{-\frac{2(q-1)}{q}} g(t)
				\bigl[\log g(t)\bigr]^{\frac{2(q-1)}{q}}
			\right) \,dt
	\]
	with $g(0)=\e^{\theta |x|_{\H}^{q}}$.
	
	By the comparison theorem, we have 
	\begin{equation}\label{Equ:Conc:62}
		\E\e^{\theta |X_{t}(x)|_{\H}^{q}} 
		\leq g(t) \leq \exp\left( C_{9}(1+t^{-\frac{q}{q-2}}) \right)
	\end{equation}
	for some constant $C_{9}>0$. 
	By inequality \eqref{MSEE:BBS:2} we have 
	\begin{equation}\label{Equ:Conc:63}
	\|P_{t}f\|_{\infty}
		=\|P_{t/2}P_{t/2}f\|_{\infty}\leq 
		C_{10} \sup_{x\in\overline{D(A)}} \E\exp
		\left[ 	\frac{C_{11}}{t^{\frac{4+r}{2+r}}} \bigl|X_{t/2}(x)\bigr|_{\H}^{\frac{2(4+r-q)}{2+r}}	\right]
	\end{equation}
	for some constants $C_{10}, C_{11}>0$. 
	
	By using Young's inequality, we see
	\begin{equation}\label{Hyperboundedness:Last:Equ}
		C_{10} \sup_{x\in\overline{D(A)}} \E\exp
		\left[ 	\frac{C_{11}}{t^{\frac{4+r}{2+r}}} \bigl|X_{t/2}(x)\bigr|_{\H}^{\frac{2(4+r-q)}{2+r}}	\right]
		\leq 
		\varepsilon \cdot |X_{t/2}^{x}|_{\H}^{q} +
		\varepsilon' t^{-\frac{q}{q-2}}	
	\end{equation}
	for arbitrary $\varepsilon>0$. By considering small enough $\varepsilon>0$, 
	it follows from the inequality \eqref{Hyperboundedness:Last:Equ} above and \eqref{Equ:Conc:62} and \eqref{Equ:Conc:63} we can obtain
	\eqref{Equ:Ultrabound:MSEE}. This proves that $P_{t}$ is ultrabounded. 
	
	Since $P_{t}$ has transition density with respect to $\mu$, we know $P_{t}$ is compact in $L^{2}(\overline{D(A)},\mu)$ for large $t>0$ 
	for these two cases ($q=2$ or $q>2$) by \cite[Lemma 3.1]{GW02}.
\end{proof}

By \cite[Section 2]{Wang09c}, 
we can get some 
heat kernel inequalities from Harnack inequality. 
Especially we have the following 
 log-Harnack inequality \eqref{log_HI:HIMSEE} and 
 an entropy-cost type inequality \eqref{Entropy:Cost:Ineq:HIMSEE}. 

\begin{theorem}\label{entropy_cost:logHI:MSEE}
Assume \hyperlink{H1:MSEE}{(H1)}--\hyperlink{H5:MSEE}{(H5)}, 
	\hyperlink{H6prime:MSEE}{(H6\,$'$)} and \eqref{ControlNormCond:MSEE}
with $q<4+r$.
Then 
\begin{enumerate}
\item For every $T>0$,  \(x,y\in \overline{D(A)}\), \(f\geq 1\), and \(f\in\BB_b(\overline{D(A)})\)
\begin{equation}\label{log_HI:HIMSEE}
	P_T(\log f)(x)\leq \log P_Tf(y) + \frac12  \Theta_{T}|x-y|^{\frac{2(4+r-q)}{2+r}}_{\H}.
\end{equation}
\item For every  $T>0$,  $f\in \BB_b^+(\overline{D(A)})$ such that $\mu(f)=1$,  
\begin{equation}\label{Entropy:Cost:Ineq:HIMSEE}
	\mu (P_T^* f \log(P_T^* f) )\leq  \frac{\Theta_{T}}{2} 
	 \inf_{\pi\in\mathscr C(\mu,f\mu)} \int_{\overline{D(A)}}\int_{\overline{D(A)}}
				 |x-y|_\H^{\frac{2(4+r-q)}{2+r}} \,d\pi,
\end{equation}
where $\mathscr C(\mu,f\mu)$ denotes all the coupling of $\mu$ and $f\mu$.
\end{enumerate}
\end{theorem}

\begin{proof}
    Inequality \eqref{log_HI:HIMSEE} follows as \cite[Proposition 2.2]{Wang09c}.
    Integrating both sides of \eqref{log_HI:HIMSEE} with respect to
$\pi\in\mathscr C(\mu,f\mu)$, and minimizing in $\pi$, we get
$$\mu((P_t^* f)\log P_t^* f) \le \mu(\log (P_tP_t^*f)) + \frac
{\Theta_T}{2} W_2(\mu, f\mu)^2.$$ 
Now we obtain \eqref{Entropy:Cost:Ineq:HIMSEE} by noting 
$$ \mu(\log(P_tP_t^*f))\le \log \mu(P_t P_t^*f)=\log
1=0$$
due to the fact that $\mu$ is invariant for $P_t$ and $P_t^*$.
\end{proof}

\section*{Acknowledgment}
This work was supported in part by the DFG through the 
International Graduate College ``Stochastics and Real World Models".
The author would like to thank the encouragement of Professors Michael R\"ockner and Feng-Yu Wang. 


\def\cprime{$'$}

\end{document}